\documentclass[12pt]{article}

\usepackage{amsmath,bm}
\usepackage{amsmath}
\usepackage{amsbsy}
\usepackage{amsbsy}
\usepackage{amssymb}
\usepackage{graphicx}

\graphicspath{{Files/}}
\usepackage{subfigure}
\usepackage{textcomp}   
\usepackage{verbatim}   
\usepackage{mathtools}
\usepackage[a4paper, total={6in, 8in}]{geometry}
\usepackage{color,soul}
\sethlcolor{cyan}
\usepackage{xspace}
\newcommand{\bz }{B\'ezier \xspace}
\newcommand{\hlb}{\textcolor{black}}
\newcommand{\hcb}{\textcolor{black}}
\newcommand{\mathbff}{\textcolor{black}}
\newcommand{\boldsymboll}{\textcolor{black}}

\fontencoding{T1}
  \fontfamily{garamond}
  \fontseries{m}
  \fontshape{it}
  \fontsize{12}{15}
\date{}
\begin{document}

\title{\hlb{$C^{1}$ Triangular Isogeometric Analysis of the von Karman Equations}}
\author{Mehrdad Zareh, Xiaoping Qian   \\
	University of Wisconsin-Madison  \\
	}


\maketitle

\begin*{\textbf{\textit{Abstract-}}}
In this paper, we report the use of rational Triangular B\'ezier Splines (rTBS) to numerically solve the von Karman equations, a system of fourth order PDEs. $C^{1}$  smoothness of the mesh, generated by triangular B\'ezier elements, enables us to directly solve the von Karman systems of equations without using mixed formulation. \textcolor{black}{Numerical results of benchmark problems show high accuracy and optimal convergence rate in $L^{1}, H^{1}$ and $H^{2}$ norm for quadratic and cubic triangular B\'ezier elements. Results of this study show that triangular isogeometric analysis can efficiently and accurately solve systems of high order PDEs.}
\end*{}
\\

\begin*{\textit{keywords:}}
Rational triangular B\'ezier splines; Isogeometric analysis; Von Karman; Optimal convergence rate; High order PDEs; Smooth mesh; Triangular elements.\footnote{This manuscript is accepted for
	publication in the Springer INdAM volume entitled "Geometric Challenges
	in Isogeometric Analysis".}
\end*{}



\section{INTRODUCTION}

Linear and nonlinear structural analysis of plates have been addressed by various formulations. In linear analysis, the most appealing theories are Kirchhoff-Love for thin plates in which transverse shear deformation is neglected, and Reissner-Mindlin theory for moderately thick plates. 

Von Karman, in 1910, introduced a system of fourth order elliptic equations to mathematically model the nonlinear behavior of plates in large deflections. \hlb{Such 4th order partial differential equations (PDEs) have also been applied in other nonlinear problems such as multiphysics modeling, ionic polymer metal composites \cite{Michopoulos2018} and  the growth of biological tissues \cite{Nelson2013}. 
Fourth order PDE and nonlinearity cause complexity in numerically solving von Karman equations \cite{Brenner2017}. Finite element method (FEM) is a widely used numerical simulation method in solid and fluid mechanics mostly because of its well-established mathematics, ability to approximate model geometries  and its generality in numerical solution, i.e. solution and its derivatives can be estimated at any location; therefore, FEM has been an attractive numerical technique to approximate the solution of von Karman equations \cite{Mallik2016, Brenner2017,Carstensen2019,Miyoshi1976,Reinhart1982}}.

For such fourth order PDEs, $C^{1}$ smoothness over the mesh is essential, which can be achieved by using relatively sophisticated finite elements (FE) such as Bogner-Fox-Schmit and Argyris FE. However, these elements are complex to implement and computationally expensive \cite{Mallik2016, Brenner2017}. \hcb{Hence, alternative methods such as smooth splines in the context of isogeometric analysis (IGA) have been explored. }

\hcb{Several techniques in the framework of IGA have been developed and used to solve high order PDEs, particularly fourth order PDEs. In \cite{Kiendl2009,Kiendl2010,Evans2020}, nonuniform rational B-Splines (NURBS) was implemented to represent shell geometries and to discretize high order PDEs of Kirchhoff–Love theory. Despite its strengths, the tensor product nature of NURBS makes them unattractive for representing complex geometries and for local mesh refinement. This problem motivated research on using multiple patches instead of single patch. Although the technique can be effective \cite{YangHS,HerremaAU,ADAM2020113403}, enforcing $C^{1}$ over the mesh constructed by multiple patches is challenging. This issue has been addressed by developing novel formulations to create $C^{1}$ splines on unstructured quadrilateral meshes \cite{May2016,Kapl2020,2018arXiv181209088K,TOSHNIWAL2017411}. Besides unstructured quadrilateral elements, triangular elements have been explored in the context of IGA. \cite{Sabin2006,Leuven2008,May2016} employed Powell-Sabin B-splines to solve the equations of fourth order PDEs and Kirchhoff–Love plate theory. Flexibility in meshing makes Powell-Sabin B-splines desirable compared with NURBS. In the same context, rational \bz triangles were used in \cite{Liu2018} to solving Kirchhoff plate problem. They impose continuity constraints by the method of Lagrange multipliers. In addition to the discussed methods, a few other strategies provide alternative tools to cope with $C^{1}$ requirement.}

\hcb{Nonconforming FE with $C^{0}$ penalty method \cite{Brenner2017} and discontinuous Galerkin FE \cite{Carstensen2019} have been employed as alternatives to eliminate the $C^{1}$-continuity requirement. Mixed finite element has also been used to solve high order PDEs \cite{Miyoshi1976,Michopoulos2018,Chen2020,Reinhart1982}; in this method additional variables are introduced into the original problem to lower the order of PDEs and, subsequently, relax the need for $C^{1}$ mesh. However, mixed FEM significantly increases the DOFs. This can cause difficulties in convergence for nonlinear problems. Moreover, introducing a new variable can change the original problem, leading to physically irrelevant solutions in some cases, e.g. geometries with re-entrant corners. In such geometries numerical solutions from mixed FEM and $C^{1}$-FEM (direct solution), for example, for a biharmonic PDE are different. In a mathematical explanation, the direct solution from $C^{1}$ FEM is in $H^2$ space; however, the solution from mixed FEM is not necessarily in $H^2$ space. This inconsistency implies that mixed FEM solution can not consistently converge to the correct solution that  is otherwise achievable with the $C^{1}$-smooth modeling of the original problem \cite{Nazarov2007,Gerasimov2012a}.} 

In this study, we employ $C^{1}$ rational Triangular B\'ezier Splines (rTBS) in the framework of triangular isogeometric analysis (TIGA) to solve von Karman equations without introducing new variables. TIGA was developed in \cite{Jaxon2014a,XIA2015292,Speleers2012,Xia2017} and has been implemented to solve various low and high order PDEs \cite{WANG2018585,Zhang2019,Speleers2012} such as Kirchhof-Love plate and shell equation \cite{Zareh2018,Zareh2019}. The availability of $C^{1}$ continuous mesh in TIGA enables us to efficiently and accurately solve the high order PDEs. This paper presents numerical examples to show the validity and optimal convergence rate in solution of von Karman equations using TIGA. 
  
\section{Rational Triangular B\'ezier Splines}
In this section, we describe some background technologies on how to achieve $C^{1}$ smooth Bezier elements. More details can be found in \cite{XIA2015292}.

We use triangular B\'ezier elements to discretize both geometry and the solution field. A B\'ezier curve is defined by Bernstein basis functions; d-degree Bernstein polynomial is given by

\begin{equation}
\label{Bernstein}
\psi_{\mathbff{I},d}(\xi) =  \frac{d!}{i!j!} \zeta^i (1-\zeta)^j, \quad \mid\mathbff{I}\mid = i + j = d.
\end{equation}
 
In this study we implement B\'ezier triangles based on bivariate Bernstein polynomials. Bivariate form of equation (~\ref{Bernstein}) describes the bivariate Bernstein polynomials:

\begin{equation}
\psi_{\mathbff{I},d}(\boldsymboll\zeta) = \frac{d!}{i!j!k!} \zeta_1^i \zeta_2^j \zeta_3^k,  \quad \mid\mathbff{I}\mid = i + j+k = d,
\end{equation}
 
\begin{figure}
\centering
\includegraphics[scale=0.4]{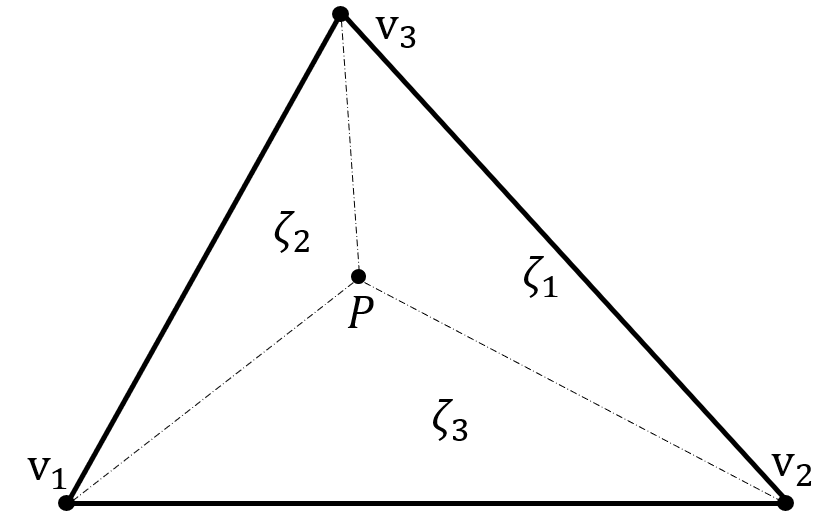}
\caption {Barycentric coordinates of a point, $P(s,t)$, in a triangle, $\tau$. $\zeta_{1}=\frac{area(V_3V_2P)}{area(V_1V_2V_3)}$, $\zeta_{2}=\frac{area(V_3V_1P)}{area(V_1V_2V_3)}$, $\zeta_{3}=\frac{area(V_1V_2P)}{area(V_1V_2V_3)}$.} 
\label{triangle} 
\end{figure}

\noindent \textbf{i} refers to a triple index i, j, k. $\zeta_{1},\zeta_{2},\zeta_{3}$ are the barycentric coordinates of a point $(s,t) \in \mathbb{R}^2$. Every points in a fixed triangle $\tau$ defined by vertices $\mathbff{v}_1,\mathbff{v}_2,\mathbff{v}_3$ (see Figure \ref{triangle}) can be uniquely defined by
\begin{equation}
(s,t) = \zeta_1\mathbff{v}_1 + \zeta_2 \mathbff{v}_2 + \zeta_3\mathbff{v}_3,\quad \zeta_1 + \zeta_2 + \zeta_3 = 1.
\end{equation} 

Now, a triangular B\'ezier patch can be defined by

\begin{equation}
\label{equ:x_xi}
\mathbff{x}(\boldsymboll\zeta) =\displaystyle\sum_{\mid \mathbff{I}\mid=d} \mathbff{p}_{\mathbff{I}} \psi_{\mathbff{I},d}(\boldsymboll\zeta),
\end{equation}
 
\noindent $\mathbff{p_{i}}$ represents the control points. By introducing the weights into the above formula, a rational B\'ezier triangle is defined by

\begin{equation}
\label{equ:x_xi2}
\mathbff{x}(\boldsymboll\zeta) =\displaystyle\sum_{\mid \mathbff{I}\mid=d} \mathbff{p}_{\mathbff{I}} \Psi_{\mathbff{I},d}(\boldsymboll\zeta),
\end{equation}
 
\noindent where
\begin{equation}
\Psi_{\mathbff{I},d}=\frac{w_{\mathbff{I}}\psi_{\mathbff{I},d}}{\displaystyle\sum_{\mid \mathbff{I}\mid=d} w_{\mathbff{I}}\psi_{\mathbff{I},d}}
\label{phi}
\end{equation}

\noindent $ w_{\mathbff{I}}$  represents the weight of the control point $\mathbff{p_{i}}$. Following the isoparametric concept, same bivariate Bernstein basis on a triangle $\tau$ with vertices $\mathbff{v}_1,\mathbff{v}_2,\mathbff{v}_3$ is used for defining a d-degree polynomial function f over $\tau$ as

\begin{equation}
f(\boldsymboll{\zeta}) = \displaystyle\sum_{\mid \mathbff{I}\mid=d} b_{\mathbff{I}} \Psi_{\mathbff{I},d}(\boldsymboll\zeta).
\end{equation}
 
\noindent The $b_{\mathbff{I}}$ (or ${b}_{ijk}$) refer to the B\'ezier ordinates of f; their corresponding array of domain points are given by
\begin{equation}
\label{equ:parnodes}
 q_{ijk} = \frac{i\mathbff{v}_1 +j\mathbff{v}_2 + k\mathbff{v}_3}{d}, \; i+j+k = d.
\end{equation}
 \\
 
\begin{figure}
\centering
\subfigure[Triangular B\'ezier patch.]{\includegraphics[scale=0.35]{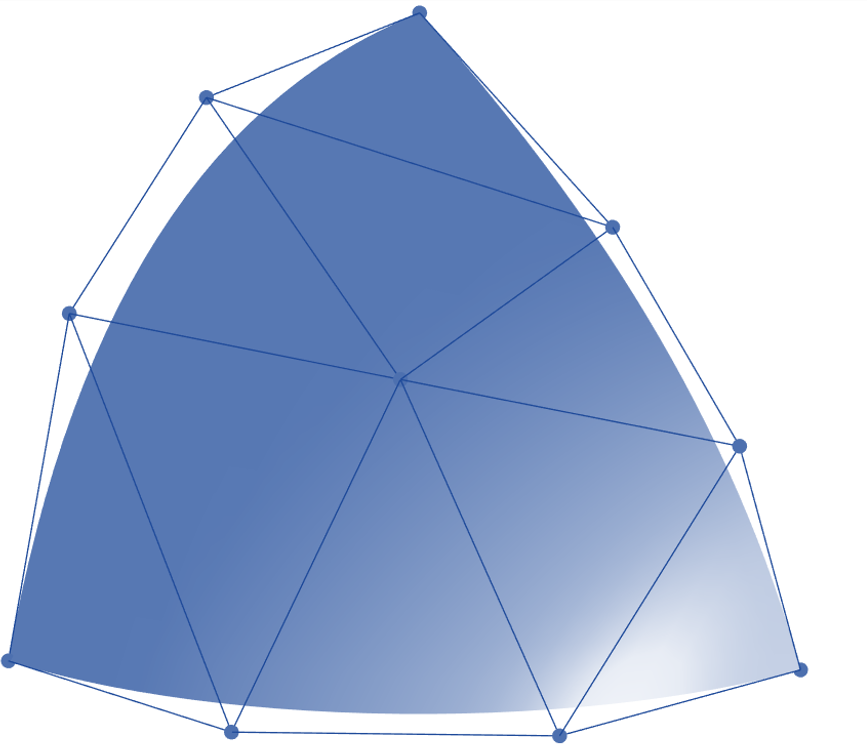}}\subfigure[Associated domain points of the B\'ezier ordinates  $b_{ijk}$ in $\left\lbrace\mathbff{v_1},\mathbff{v_2},\mathbff{v_3}\right\rbrace$.   ]{\includegraphics[scale=0.4]{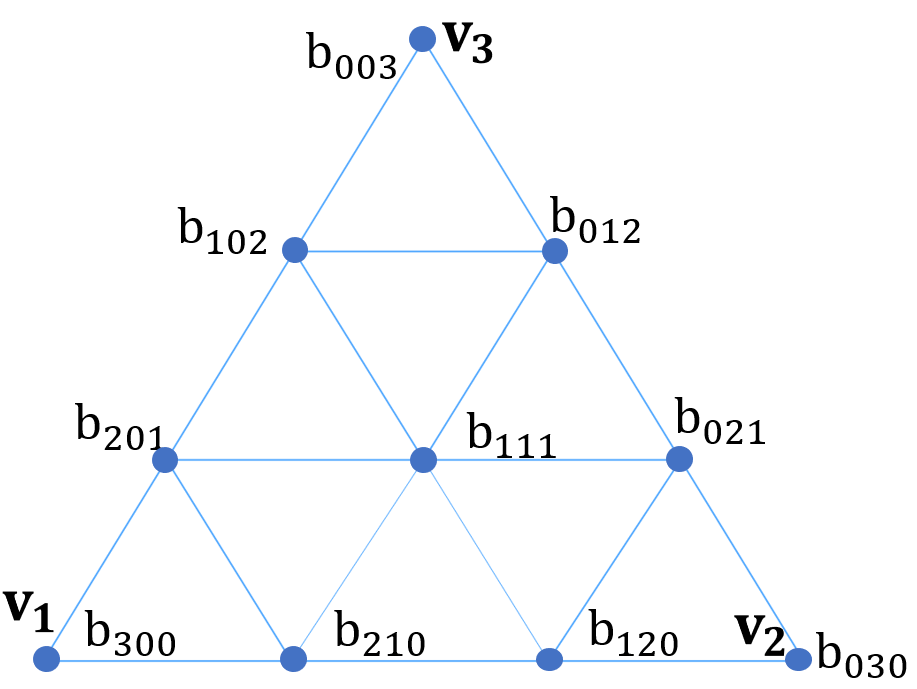}}
\caption{Triangular B\'ezier patch and domain points.}
\label{bpatch} 
\end{figure}

The control polygon of the function f is defined by the points $(q_{ijk}, b_{ijk})$. Figure \ref{bpatch}
presents an example of a triangular B\'ezier patch and the corresponding domain points of the B\'ezier ordinates.

Due to the presence of second order PDEs in weak form of von Karman equations, continuity within and between elements is the key requirement in applying finite element method for this equation. One way to relax this requirement is introduction of an additional variable into the equation, which can cause expensive computation and inaccuracy in some cases. In this study, we use $C^{1}$ triangular isogeometric method; therefore, we can achieve the required continuity. We describe how high continuity is attained between triangular B\'ezier patches. Two degree-d polynomials $f$ and $\tilde{f}$ join r times differentiably across the interface of two triangles $\tau = \left\lbrace \mathbff{v}_1,\mathbff{v}_2,\mathbff{v}_3 \right\rbrace$  and ˜ $\tilde{\tau} = \left\lbrace \mathbff{v}_4,\mathbff{v}_2,\mathbff{v}_3 \right\rbrace$ if and only if \cite{lai_schumaker_2007};$ \left( j+k+\rho = d,\ \rho=0,...,r,\right) $

\begin{equation}
\label{equ:parnodes2}
 \tilde{b_{\rho,j,k }}- \sum \frac{\rho!}{\mu!\nu!\kappa!}b_{\mu,k+\nu,j+\kappa}\zeta_{1}^\mu,\zeta_{2}^\nu,\zeta_{3}^\kappa = 0, 
\end{equation}

\begin{figure}
\centering
\includegraphics[scale=0.5]{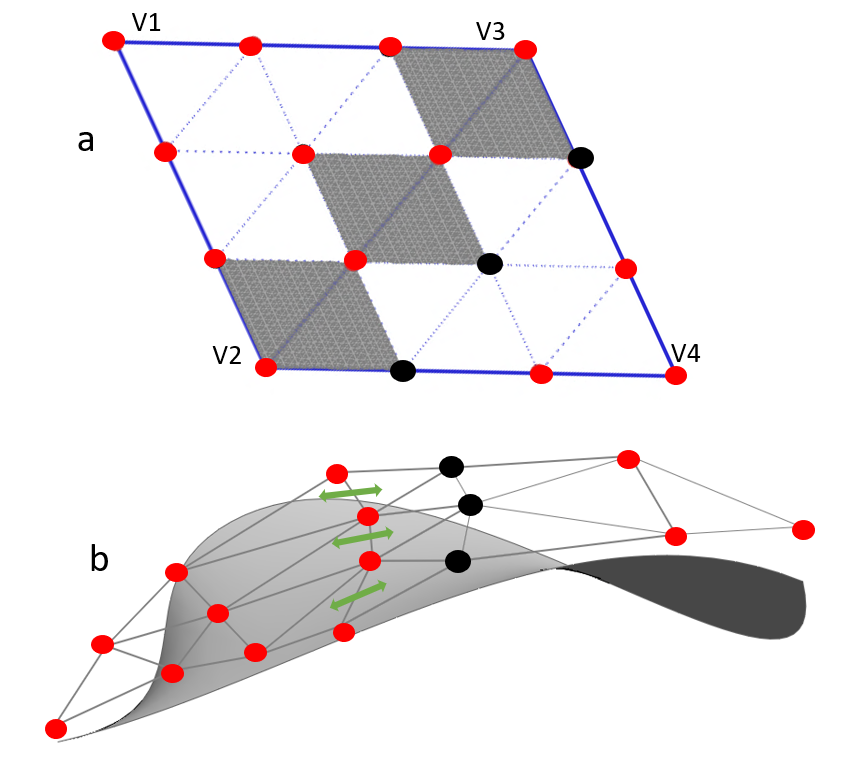}
\caption {\hcb{Dependent (black) and independent nodes (red). (a) Two domain triangles, continuity constraints are imposed on the shaded triangles. (b) two B\'ezier patches with $C^{1}$ continuity, the control points in the triangle pairs (indicated by $\longleftrightarrow)$ are coplanar. For clarity in visualization, the control net is shifted up.}} 
\label{cont} 
\end{figure}

\noindent $\zeta_{1},\zeta_{2},\zeta_{3}$ represent the barycentric coordinates of vertex v4 in relation to triangle $\tau$. Figure \ref{cont} shows two triangular B\'ezier patches with $C^{1}$ continuity across the boundary of patches. The red solids are free nodes of which values are independently computed; value of black solids (dependent nodes) are determined by applying the continuity constraints between red and blacked nodes . The continuity constraints are applied over the gray area; moreover, this figure demonstrates the coplanarity of the control points in these triangle pairs.
 
\begin{figure}
\centering
\subfigure[Triangulated physical domain]{\includegraphics[scale=0.33]{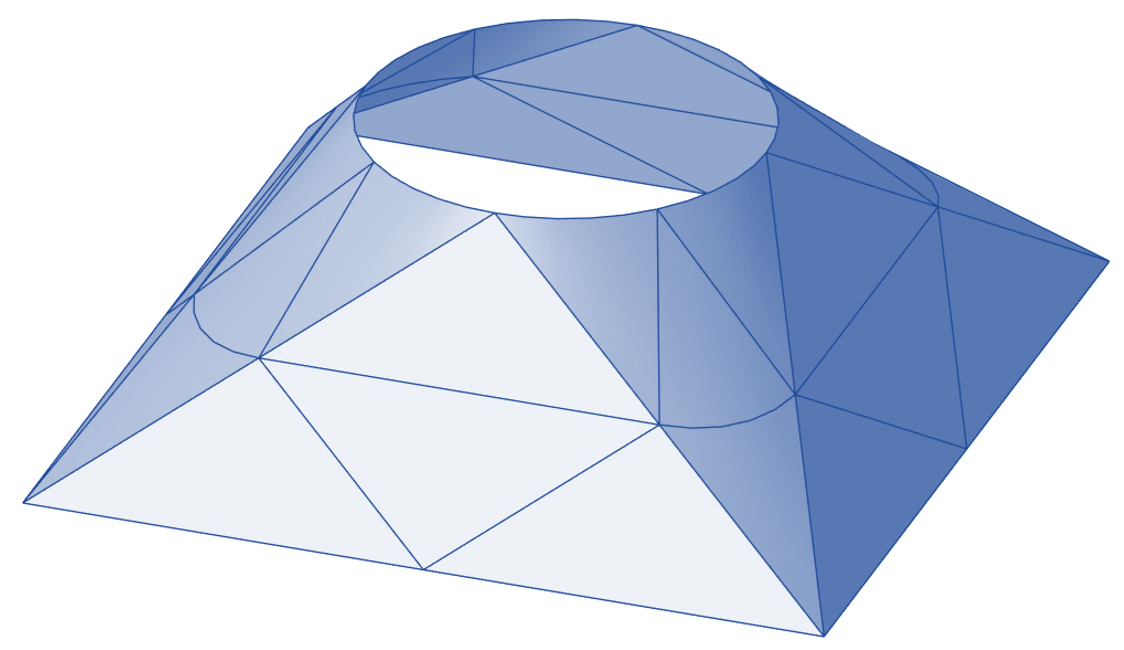}}
\subfigure[Triangulated parametric domain]{\includegraphics[scale=0.3]{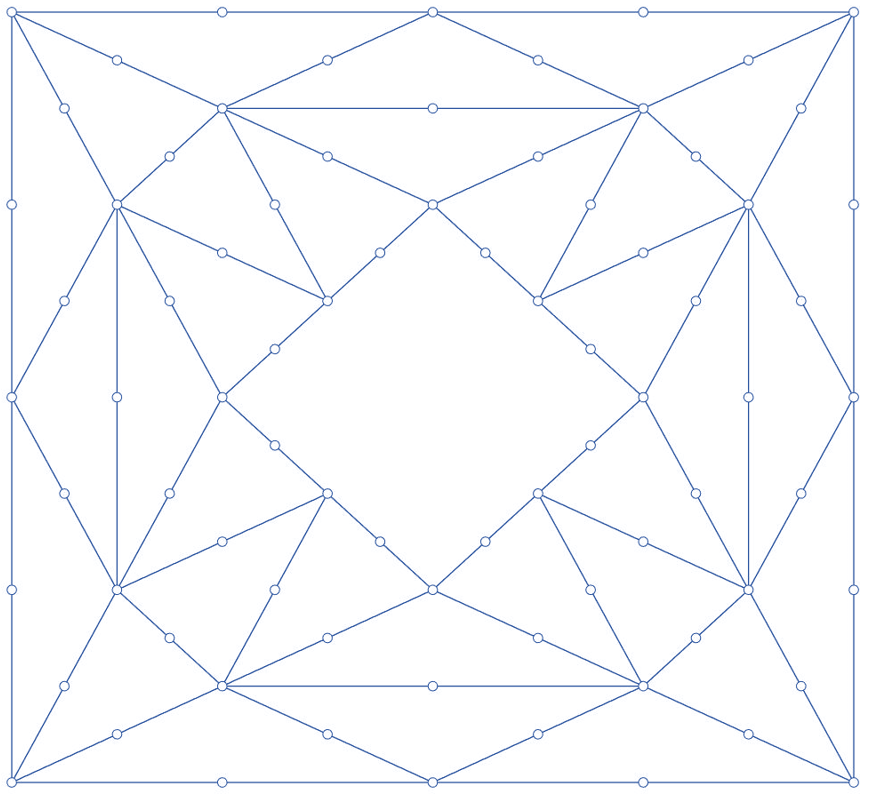}}
\caption{Illustration of a Triangulated physical and parametric Domain}
\label{param} 
\end{figure}

 Having a parametric domain $\hat{\Omega }$ and its triangulation $\hat{T}$ (Figure \ref{param}), the spline spaces of piecewise d-degree polynomials $\hat{T}$ are defined by \cite{lai_schumaker_2007};

\begin{equation}
\label{equ:parnodes3}
\mathbff{S}_{d}^{r}(\hat{T})=\left\lbrace f \in \mathbff{C}^{r}(\hat{\Omega}):f|_{\tau} \in \mathbff{P} \ \forall \tau \in \hat{T}\right\rbrace , 
\end{equation}

\noindent $\tau$ is an arbitrary triangle in $\hat{T}$, r and $\mathbff{P}$ are the continuity order of the spline over $\hat{\Omega }$ and space of polynomials of degree d. When spline has higher smoothness across some edges or at some vertices, it is called superspline; the associated space is given by~\cite{lai_schumaker_2007};

\begin{equation}
\label{equ:parnodes4}
\mathbff{S}_{d}^{r, \rho}(\hat{T})=\left\lbrace f \in \mathbff{S}_{d}^{r}(\hat{T}):f \in \mathbff{C}^{r_{v}}(V) \ \forall v \in V \ \& \ f \in \mathbff{C}^{r_{e}}(e)\ \forall e \in E\right\rbrace.  
\end{equation}

All vertices and edges are represented by V and E in $\hat{T}$ and $\rho:= \left\lbrace \rho_{v} \right\rbrace_{v \in V } \cup \left\lbrace \rho_{e}\right\rbrace_{e \in E } $with $ r \le \rho_{v}, \rho_{e}\le d$ for each $v \in V$ and $e \in E$. 

\hcb{Before applying continuity constraint, a minimal determining set (MDS) is framed, such set contains all free domain points. In this study we use direct construction (DC) method to construct MDS, in which a set of free domain points are directly chosen based on the connectivity of the triangle elements. \cite{XIA2015292} details the DC method and the alternative method of Gaussian elimination to build MDS.}
	
\begin{figure}
\centering
\subfigure[PS macro-element]{\includegraphics[width=0.4\textwidth]{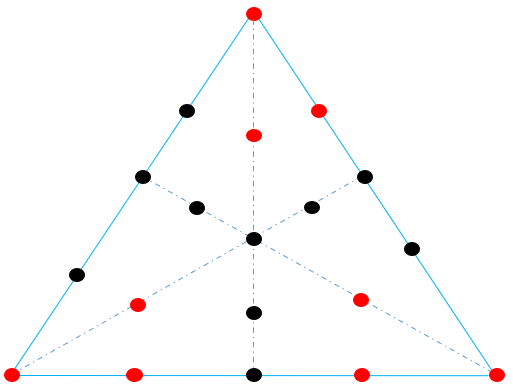}}\qquad
 \subfigure[CT macro-element.]{\includegraphics[width=0.4\textwidth]{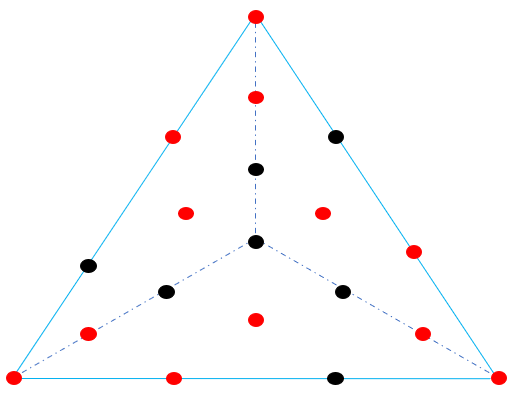}}
 
\caption{Splitting methods used in this study. red and black nodes represent independent and dependent nodes, respectively.}
\label{CTPS}
\end{figure}

\hcb{Imposing condition \ref{equ:parnodes2} directly on the triangles is a conventional technique to create $C^{r}$ spline spaces on a triangulated domain $\hat{\Omega}(\hat{T}) $. Despite being straightforward, this direct method requires the degree of the polynomial to be much higher than r, i.e.  $d \ge  3r+2$ \cite{XIA2015292}. In this work, triangles in $\hat{T}$ are split into multiple microtriangles before imposing the continuity constraints on the microtriangles. Clough-Tocher (CT) and Powell-Sabin (PS) methods are implemented for splitting; these methods do not need high-degree polynomials to provide continuity, e.g. CT cubic elements and PS quadratic elements can create $C^{1}$ mesh. In the CT splitting method, each vertex of a triangle is connected to its centroid point, forming three micro-triangles. The PS method splits each macro-triangle into six micro-triangles with centroid point as the interior split point. Edges are then bisected (see Figure \ref{CTPS}).}

\hcb{In order to evaluate convergence rate, we use smooth-refine-smooth procedure \cite{Xia2017}. The key point of this method is the relation between the sequence of triangulation in refinements, such that the refined triangulation (mesh) is the subset of the pre-refined triangulation. In this technique, continuity constraint is applied before refining the mesh in order that a smooth pre-refinement map is constructed; then, further refinement followed by $C^{r}$ continuity constraints are applied. The smooth-refine-smooth procedure ensures that control points do not relocate as they have already satisfied the continuity conditions. The obtained mesh is $C^r$ smooth, and the geometric map remains the same; therefore, inconsistency is avoided for all subsequent refinements.}

\section{Formulations and Methods}
\hcb{This section presents the von Karman equations, the corresponding mixed formulation, the weak formulations and finite element discretization. As discussed in Introduction, $C^{1}$ mesh is required for directly solving original von Karman equations due to the high order PDEs. On the other hand, such high order derivatives do not appear in the weak form of mixed FEM; therefore,  $C^{0}$ mesh, available in conventional FEM, is sufficient for solving mixed formulation.}

\subsection{Von Karman equations}

The canonical form of von Karman equations with simply supported boundary condition is described by

\begin{equation}
\label{vonk}
\begin{cases}
{\Delta^{2} u=[u, v]+f}, \\
{\Delta^{2} v=-\frac{1}{2}[u, u]},
\end{cases}
\end{equation}

\noindent where bi-harmoinc operator and the von Karman bracket are defined by\\
\\$\displaystyle\Delta^{2} u \equiv \frac{\partial^{4} u}{\partial x^{4}}+2 \frac{\partial^{4} u}{\partial x^{2} \partial y^{2}}+\frac{\partial^{4} u}{\partial y^{4}}$,
\\
$\displaystyle[u, v] \equiv \frac{\partial^{2} u}{\partial x^{2}} \frac{\partial^{2} v}{\partial y^{2}}+\frac{\partial^{2} u}{\partial y^{2}} \frac{\partial^{2} v}{\partial x^{2}}-2 \frac{\partial^{2} u}{\partial x \partial y} \frac{\partial^{2} v}{\partial x \partial y}.
$\\

The above equations have fourth order PDEs. Because of their complications in finite element formulation, mixed form of von Karman equations has been extensively used, which does not have high order PDEs. For completeness in the context, next session presents the  mixed form; however, we note that the focus of this article is on using $C^{1}$ mesh for von Karman equations with fourth order PDEs.

\subsection{Mixed form of von Karman equations}
By introducing an additional variable, $\bar{u}$, into the equation (\ref{vonk}), the mixed formulation is defined by

\begin{equation}
		\label{mixvonk}
		\begin{cases}
			{\Delta u=\bar{u}},\\
			{\Delta \bar{u}=[u, v]+f},\\
			{\Delta v=\bar{v}},\\
			{\Delta \bar{v}=-\frac{1}{2}[u, u],}
		\end{cases}
\end{equation}

\noindent where $\Delta$ is the Laplacian operator, and the von Karman bracket was defined in the previous session.

As we can observe, the above formulation does not have high order PDEs, i.e. fourth order PDEs do not appear in the equations.
Next section provides the finite element discretisation of the above equations.

\subsection{Finite element discretisation}

The strong form of von Karman equations are multiplied by test functions to obtain the weak formulations. 

\subsubsection{The weak form of mixed formulation}
 The following equation represents the weak form of (\ref{mixvonk}). (for $f\in L^2(\Omega)$ and $ u, v, \bar{u}, \bar{v}\in S:=H_0^1(\Omega) $)
\hlb{
	\begin{equation}
		\label{weakM}
		\begin{cases}
			{a(u,w)+m(\bar{u},w)=0, \forall w \in S},\\
			{a(\bar{u},\mu)+b(u,v,\mu)+b(v,u,\mu)=l(w), \forall \mu \in S},\\
			{a(v,\tau)+m(\bar{v},\tau)=0, \forall \tau \in S},\\
			{a(\bar{v},\phi)-b(u,u,\phi)=0, \forall \phi \in S},
		\end{cases}
	\end{equation}
}

\noindent\hlb{
	the bilinear $a$, trilinear $b$ and $l$ are given by ($\forall \beta, \gamma, \varphi \in S$, and $ \operatorname{cof}\left(D^{2}\beta\right)$ is the cofactor matrix of the Hessian of $\beta$)\\ 
	$\displaystyle 
	a(\beta, \gamma):=-\int_{\Omega} \bigtriangledown \beta. \bigtriangledown \gamma d \Omega.\\
	$
}
$\displaystyle b(\beta, \gamma, \varphi):=\frac{1}{2} \int_{\Omega} \operatorname{cof}\left(D^{2} \beta\right) D \gamma \cdot D \varphi d\Omega,\\$
$\displaystyle  l(w)=\int_{\Omega}f\cdot w d\Omega.\\$

\subsubsection{The weak form for the original formulation }
The weak formulation of (\ref{vonk}) is given by (for $f\in L^2(\Omega)$ and $ u, v\in S:=H_0^2(\Omega) $)

\begin{equation}
\label{weak}
\begin{cases}
{\hat{a}(u,w)+b(u,v,w)+b(v,u,w)=l(w), \forall w \in S},\\
{\hat{a}(v,\phi)-b(u,u,\phi)=0, \forall \phi \in S}.
\end{cases}
\end{equation}

The bilinear $\hat{a}$ is defined by 
\\$\displaystyle 
\hat{a}(\beta, \gamma):=\int_{\Omega} D^{2} \beta: D^{2} \gamma d \Omega,\\$

As it can be observed, weak form (\ref{weak}) has second derivative of solution. Therefore, $C^{1}$ continuity between and within elements is needed in order to obtain a compatible finite element method to meet computability. We proceed with (\ref{weak}) to describe the next steps in solving von Karman equations; the same procedure can be applied to (\ref{weakM}).  

After discretization, following element matrices are defined;
\begin{equation}
K_{e}=\int_{\Omega_{e}}H_{e}^{T}H_{e} d\Omega,
\label{stiff}
\end{equation}

\begin{equation}
	F_{e}=\int_{\Omega_{e}}f\Phi_{i}d \Omega,
\end{equation}

\begin{equation}B_{\mathbff{e}}^{\mathbff{m}}:=\left[b_{i, j}^{m}\right]_{1 \leq i, j \leq n}, 1 \leq m \leq n,
\label{bstiff}
\end{equation}

\noindent where (n=number of nodes in the element)
\[H_{e}= 
\begin{bmatrix}
    \Phi_{1,xx}       & \Phi_{2,xx} &...  & \Phi_{n,xx} \\
   \Phi_{1,yy}       & \Phi_{2,yy} &...  &\Phi_{n,yy} \\
\end{bmatrix},
\]

\begin{equation}\begin{aligned}
b_{i, j}^{k}=& \frac{1}{2} \int_{\Omega}\left(\left(\frac{\partial^{2} \Phi_{i}}{\partial y^{2}} \frac{\partial \Phi_{j}}{\partial x}-\frac{\partial^{2} \Phi_{i}}{\partial x \partial y} \frac{\partial \Phi_{j}}{\partial y}\right) \frac{\partial \Phi_{k}}{\partial x}\right.\\
&\left.+\left(\frac{\partial^{2} \Phi_{i}}{\partial x^{2}} \frac{\partial \Phi_{j}}{\partial y}-\frac{\partial^{2} \Phi_{i}}{\partial x \partial y} \frac{\partial \Phi_{j}}{\partial x}\right) \frac{\partial \Phi_{k}}{\partial y}\right) d \Omega.
\end{aligned}\end{equation}

In this study, the displacement, $u$, and airy stress, $v$, of the structure are polynomial function. Following Galerkin method and FE discretization, the weighted basis function introduced in equation (\ref{phi}) is plugged in equation (\ref{stiff},\ref{bstiff}).

The element matrices form the global matrices, e.g., $\mathbff{F}$ is assembled from $F_{e}$. Finally, The discrete form is given by

\begin{equation}\begin{cases}
{\mathbff{K} \boldsymboll{u}+\left(\boldsymboll{u}^{T} \mathbff{B}^{\mathbff{m}} \boldsymboll{v}\right)_{1 \leq \mathbff{m} \leq N}+\left(\boldsymboll{v}^{T} \mathbff{B}^{\mathbff{m}} \boldsymboll{u}\right)_{1 \leq \mathbff{m} \leq N}=\mathbff{F}}, \\
{\mathbff{K} \boldsymboll{v}-\left(\boldsymboll{u}^{T} \mathbff{B}^{\mathbff{m}} \boldsymboll{u}\right)_{1 \leq \mathbff{m} \leq N}=0},
\end{cases}\end{equation}

\noindent where

$$\begin{aligned}
\left(\boldsymboll{u}^{T} \mathbff{B}^{\mathbff{m}} \boldsymboll{v}\right)_{1 \leq \mathbff{m} \leq N}:=\left[\boldsymboll{u}^{T} \mathbff{B}^{\mathbff{1}} \boldsymboll{v}\right.&\left.\boldsymboll{u}^{T} \mathbff{B}^{2} \boldsymboll{v} \quad \cdots \quad \boldsymboll{u}^{T} \mathbff{B}^{\mathbff{N}} \boldsymboll{u}\right]^{T}, \\
\boldsymboll{u}=\left(u_{1}, \cdots, u_{N}\right)^{T}, & \boldsymboll{v}=\left(v_{1}, \cdots, v_{N}\right)^{T}.
\end{aligned}$$

Due to the nonlinearity, we use Newton’s iterative methods based on the following formulation. 

\begin{equation}\begin{array}{c}
\mathbff{R}(\boldsymboll{u}, \boldsymboll{v})=\left(\mathbff{K} \boldsymboll{u}+\left(\boldsymboll{u}^{T} \mathbff{B}^{\mathbff{m}} \boldsymboll{v}\right)_{1 \leq \mathbff{m} \leq N}+\left(\boldsymboll{v}^{T} \mathbff{B}^{\mathbff{m}} \boldsymboll{u}\right)_{1 \leq \mathbff{m} \leq N}-\mathbff{F}, \mathbff{K} \boldsymboll{v}-\left(\boldsymboll{u}^{T} \mathbff{B}^{\mathbff{m}} \boldsymboll{u}\right)_{1 \leq \mathbff{m} \leq N}\right).
\end{array}\end{equation}

The solution can be compacted into $\mathbff{Z}^{n}=(\mathbff{u}^n,\mathbff{v}^n)$; The Jacobian of $\mathbff{R}$ is needed; it is defined by ($\mathbff{Bv}$ and $\mathbff{Bu}$ are computed by using $\left(\boldsymboll{u}^{T} \mathbff{B}^{\mathbff{m}} \boldsymboll{v}\right)_{1 \leq \mathbff{m} \leq N}$ )

\[\mathbff{J}(\mathbff{Z})= 
\begin{bmatrix}
    \mathbff{K}+\mathbff{Bv} & \mathbff{Bu}\\
    -\mathbff{Bu} & \mathbff{K}\\
\end{bmatrix},
\]

 \noindent then, in each iteration:

\begin{equation}
\mathbff{Z}^{n+1}=\mathbff{Z}^{n}+\mathbff{\delta}^{n},
\end{equation}

\noindent$\mathbff{\delta}^{n}$ is obtained by solving the following system;
\begin{equation}
\mathbff{J}(\mathbff{Z}^{n})\mathbff{\delta}^{n}=-\mathbff{R}(\mathbff{Z}^{n}).
\end{equation}

\section{Numerical Results}

In this section numerical examples are presented. First, the present method is verified against analytical solution. Also, convergence rate is investigated. Finally, we compare the results obtained from both mixed FEM and isogeometric analysis. 

\subsection{Convergence Study}

\begin{figure}
\centering
\subfigure[C1 Mesh]{\includegraphics[width=0.375\textwidth]{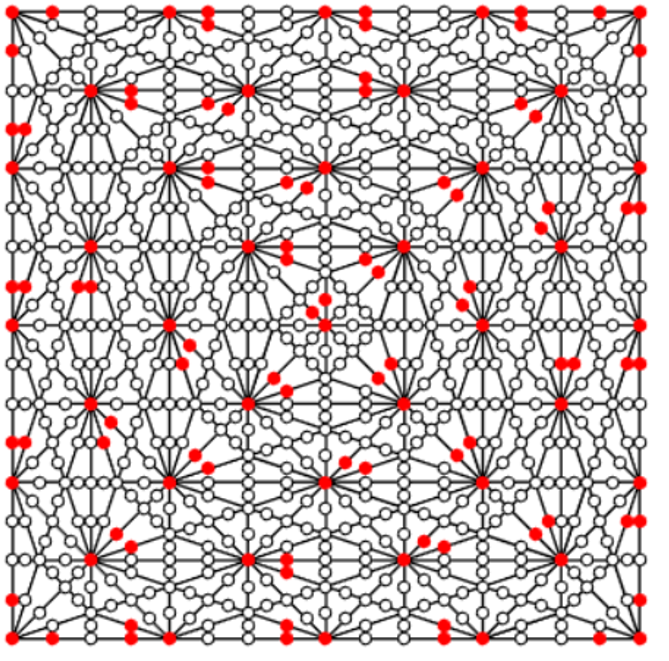}}\qquad
 \subfigure[Numerical result for $u$]{\includegraphics[width=0.5\textwidth]{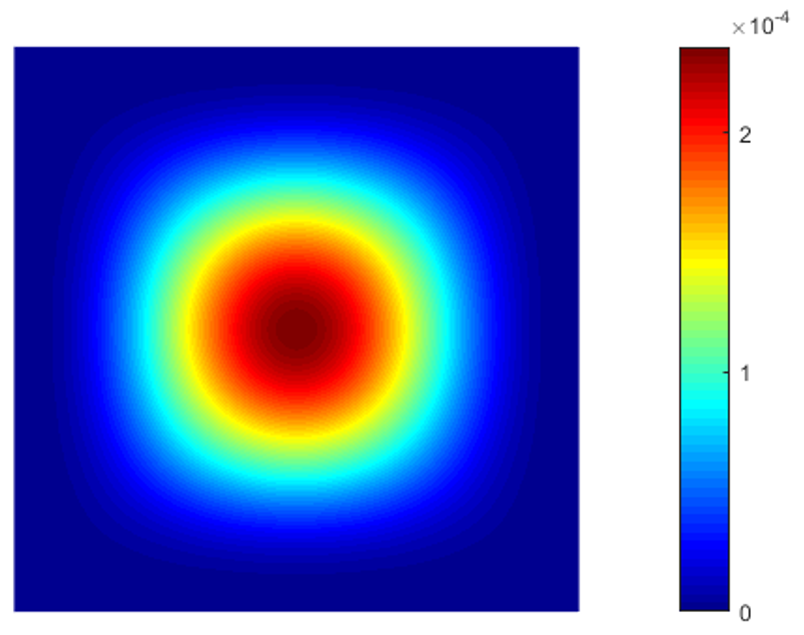}}
 
\caption{Benchmark Example: (a) mesh with 432 elements created by PS method; the red solid nodes are independent nodes and white nodes are dependent nodes. (b) numerical solution for u using $C^{1}$ mesh.}
\label{resultu} 
\end{figure}

For benchmark problem, a unit square with simply supported boundary condition is considered. The following exact solution is used in equation (\ref{vonk}) to obtain $f$ and $g$ (g is added to the rhs of the bottom of equations \ref{mixvonk} \ref{vonk}). Therefore, in numerical method, $f$ and $g$ are implemented as inputs and $u$ and $v$ are computed. \hlb{Figure \ref{resultu}  shows the mesh and numerical solution. This mesh is built from 432 elements (microtriangles) created by PS method; the red solid nodes are independent nodes (free) of which values determine the values of 
white nodes (dependent nodes) through the continuity constraints, as explained in previous section.}

\begin{equation}
\begin{aligned}
&u=x^{3}(1-x)^{3} y^{3}(1-y)^{3},\\
&v=\sin (\pi x) \sin (\pi y).
\end{aligned}
\end{equation}

We consider three error norms in convergence analysis, $L^{2}$ norm, $H^{1}$ and $H^{2}$ seminorm. First, we define the following error function:

\begin{equation}
	\label{errorf}
		e=u^{h}-u^{a},
\end{equation}

\noindent where $u^{h}$ and $u^{a}$ represent the numerical and exact solutions. We use $E$ in computing the error norms.

\begin{figure}
\centering
\includegraphics[scale=0.58]{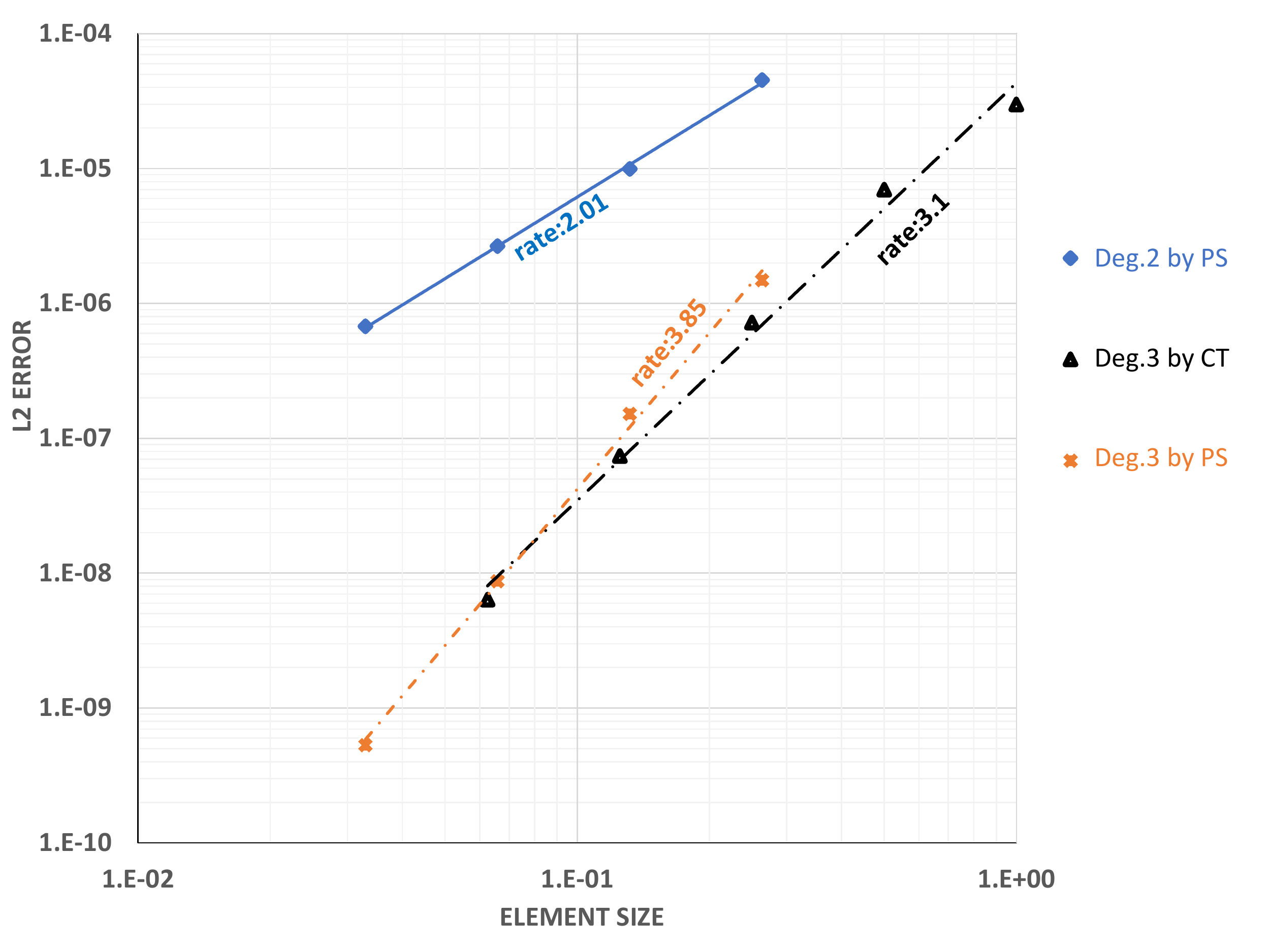}
\caption{For L2 error, the optimal convergence rate (quadratic) is obtained for TIGA using quadratic triangular elements. 
Using cubic elements, the convergence rate is 3 for TIGA-CT and 4 for TIGA-PS (optimal rate is expected to be 4)}
\label{convl2} 
\end{figure}

Convergent plot for $L^{2}$ error is illustrated in Figure \ref{convl2}; we compute this norm by
	
\begin{equation}
	\|e\|_{L^{2}(\Omega)}=\sqrt{\int_{\Omega}e\cdot ed\Omega},\\
\end{equation}
	
The optimal convergence rate is expected to be as follows: rate of 2 for quadratic elements and rate of 4 for cubic elements. We can observe in Figure \ref{convl2} that the optimal convergence rate is obtained for both quadratic (rate of 2) and cubic elements (rate of 4) created by PS method. We note that because of quadratic essence in the PS method, rate of three ($p+1=3$) is not attainable in quadratic elements \cite{Leuven2008,Bartezzaghi2015}.

\begin{figure}
\centering
\includegraphics[scale=0.58]{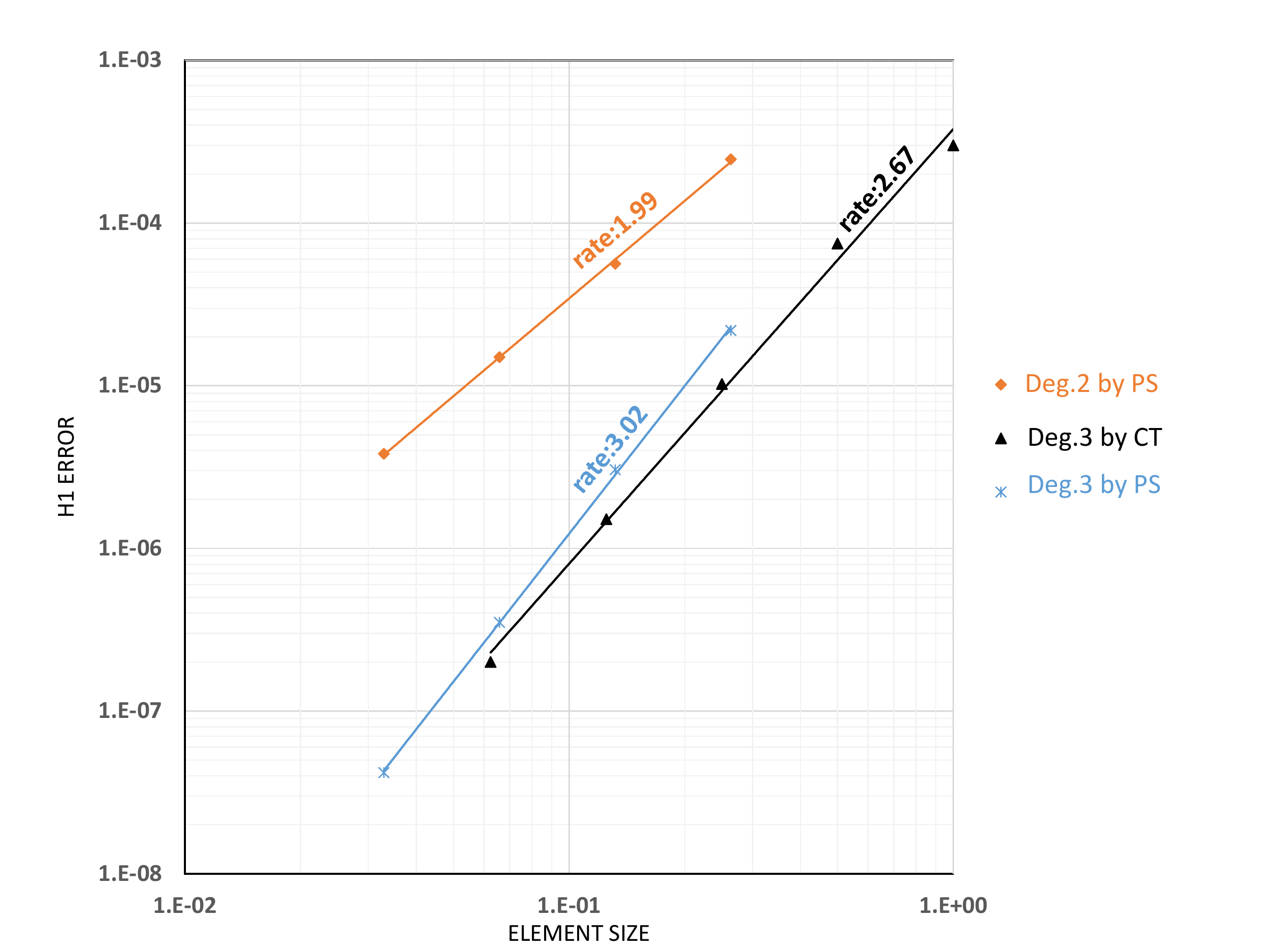}
\caption{\hcb{$H^{1}$ error analysis, the optimal convergence rates are 2 and 3 for quadratic and cubic elements, respectively. We can observe that these optimal convergence rates are obtained for quadratic elements (rate=2), and for PS cubic elements (rate=3). For cubic elements created by CT sub-optimal rate, 2.67, is obtained.}}
\label{H1rate} 
\end{figure}

Figure \ref{H1rate} presents the convergence plot for $H^{1}$ error analysis. For this analysis, we use the seminorm given by
	
\begin{equation}
	\mid e\mid_{H^{1}(\Omega)}=\sqrt{\int_{\Omega}(\nabla e\cdot \nabla e) d\Omega}.\\
\end{equation}
	
The optimal convergence rate is expected to be 2 and 3 for quadratic and cubic elements, respectively. We can observe in Figure \ref{H1rate} that the optimal convergence rate is obtained for quadratic elements (rate of 2); for cubic elements, sub-optimal rate, 2.67, is obtained when triangular elements are created by CT method; however, for cubic elements created by PS the optimal rate, 3, is obtained.

\begin{figure}
	\centering
	\includegraphics[scale=0.58]{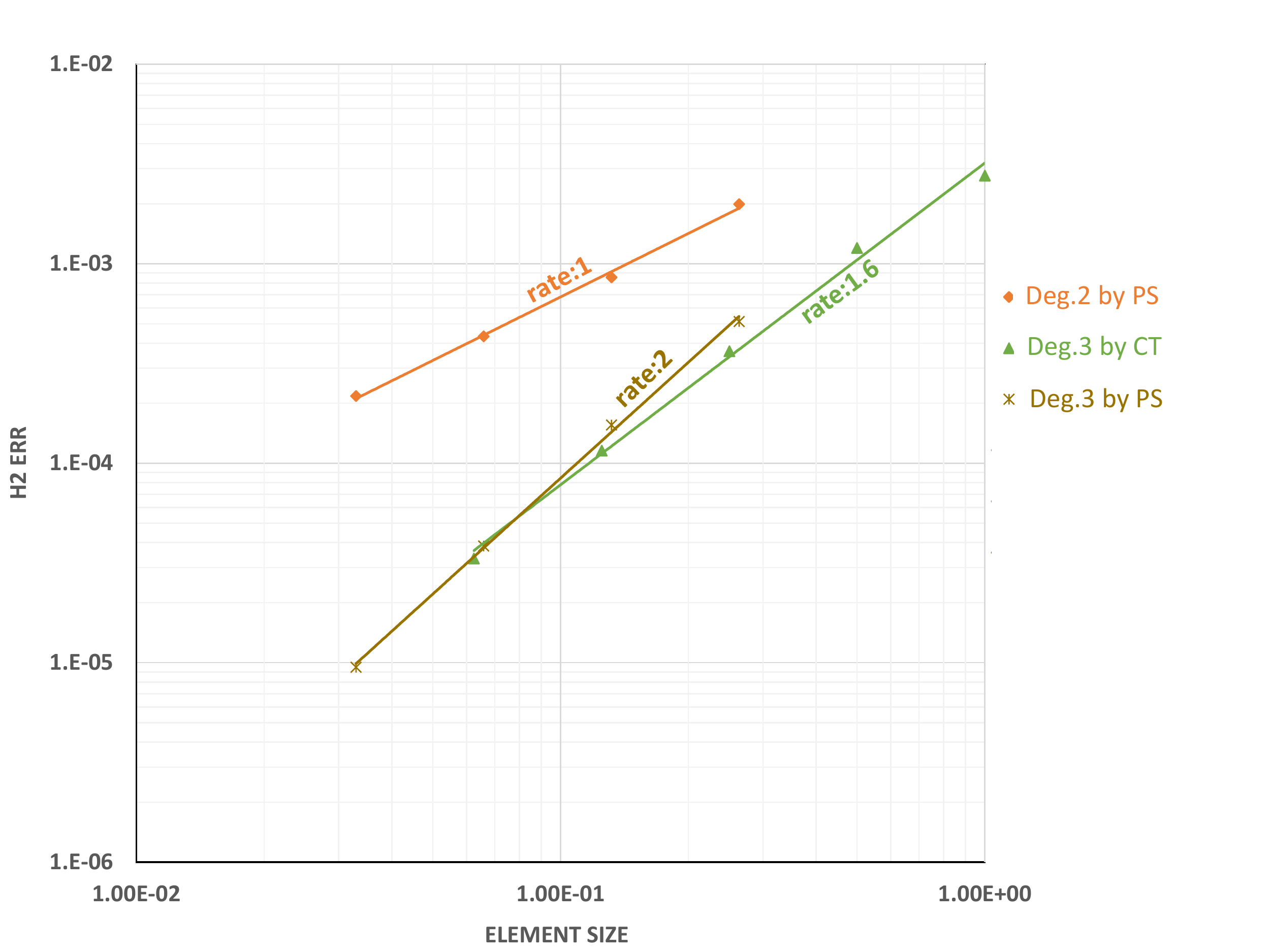}
	\caption{\hcb{$H^{2}$ error analysis, the optimal convergence rates are 1 and 2 for quadratic and cubic elements, respectively. We can observe that these optimal convergence rates are obtained for quadratic elements (rate=1), and PS cubic elements (rate=2). For cubic elements created by CT sub-optimal rate, 1.6, is obtained.}}
	\label{H2rate} 
\end{figure}

\hcb{Convergence plot for $H^{2}$ error analysis is illustrated in Figure \ref{H2rate}. For this analysis, we use the seminorm defined by}

\begin{equation}
	\mid e\mid_{H^{2}(\Omega)}=\sqrt{\int_{\Omega}(d^2 e\cdot d^2 e) d\Omega},\\
\end{equation}

\noindent where $d^2 e$ is given by
$  \displaystyle  d^2 e=\begin{bmatrix} 
\smallskip \frac{\partial^{2} e}{\partial x^{2}} 
\\ \frac{\partial^{2} e}{\partial x^{2}}\end{bmatrix}$.\\

\hcb{The optimal convergence rate is expected to be 1 and 2 for quadratic and cubic elements, respectively. We can observe in Figure \ref{H2rate} that the optimal convergence rate is obtained for quadratic elements (rate of 1); for cubic elements, sub-optimal rate, 1.6, is obtained when triangular elements are created by CT method; however, the optimal rate, 2, is obtained for cubic elements created by PS.}


\subsection{Efficiency}
\begin{figure}
\centering
\includegraphics[scale=0.6]{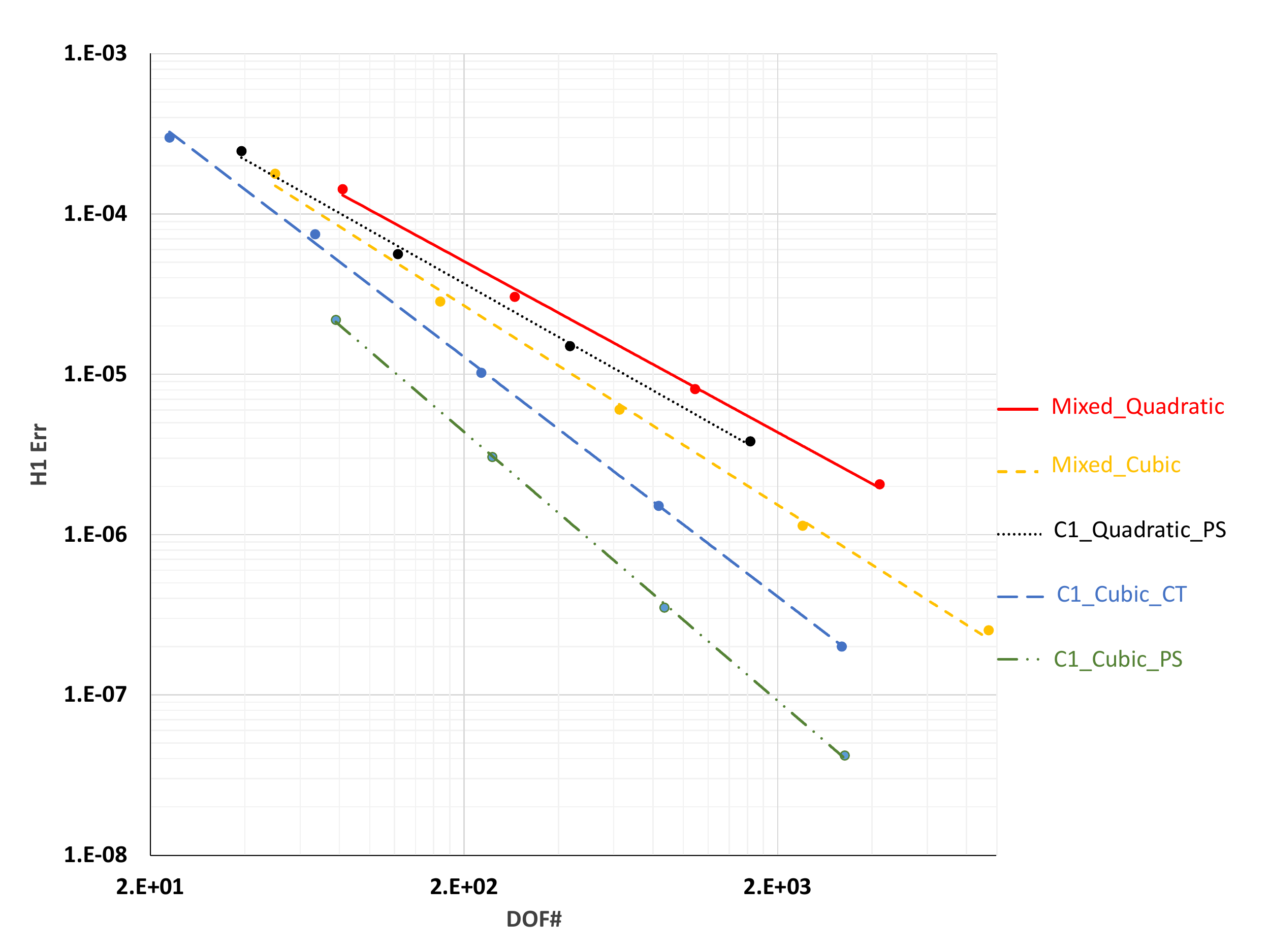}
\caption{H1 error analysis with respect to number of DOF; Assuming H1 error is expected to be $10^{-6}$ , for quadratic elements, Mixed FEM needs (7688 DOFs) ~48\% more DOFs than TIGA does (5204 DOFs). For cubic elements, Mixed FEM needs (2600 DOFs) 150\% and 490\% more DOFs than TIGA-CT (1040 DOFs) and TIGA-PS (440 DOFs) does, respectively.}
\label{H1DOF} 
\end{figure}

Figure \ref{H1DOF} compares the efficiency of $C^{1}$ TIGA with  mixed FEM; one can observe that $C^{1}$ TIGA is much more efficient than mixed FEM, e.g., assuming H1 error is expected to be $10^{-6}$ , for quadratic elements mixed FEM needs (7688 DOFs) ~48\% more DOFs than TIGA does (5204 DOFs). For cubic elements, Mixed FEM needs (2600 DOFs) 150\% and 490\% more DOFs than TIGA-CT (1040 DOFs) and TIGA-PS (440 DOFs) does, respectively. 

\subsection{Accuracy and reliability in complex geometries}
\begin{figure}
\centering
\subfigure[$C^{0}$ Mesh]{\includegraphics[scale=0.35]{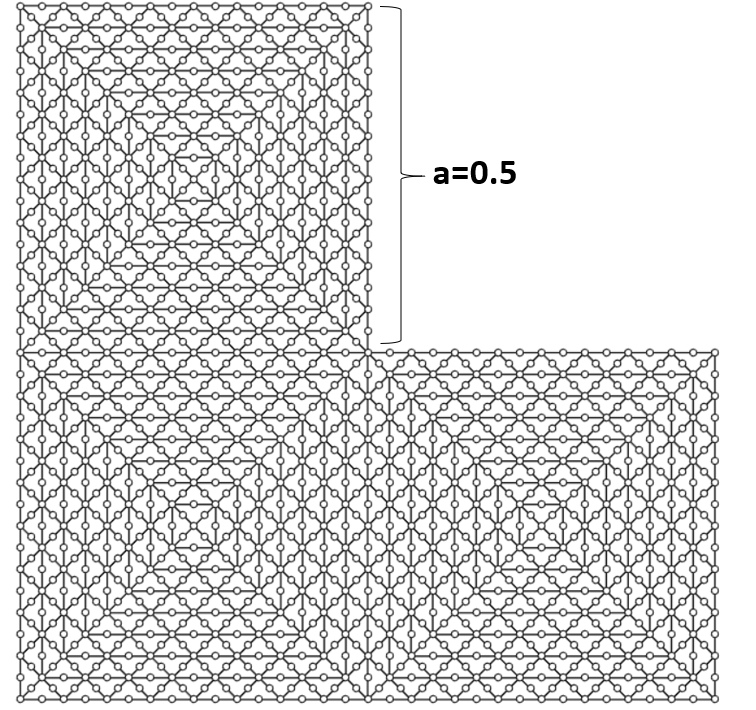}}\hspace{5mm}%
\subfigure[$C^{1}$ Mesh]{\includegraphics[scale=0.35]{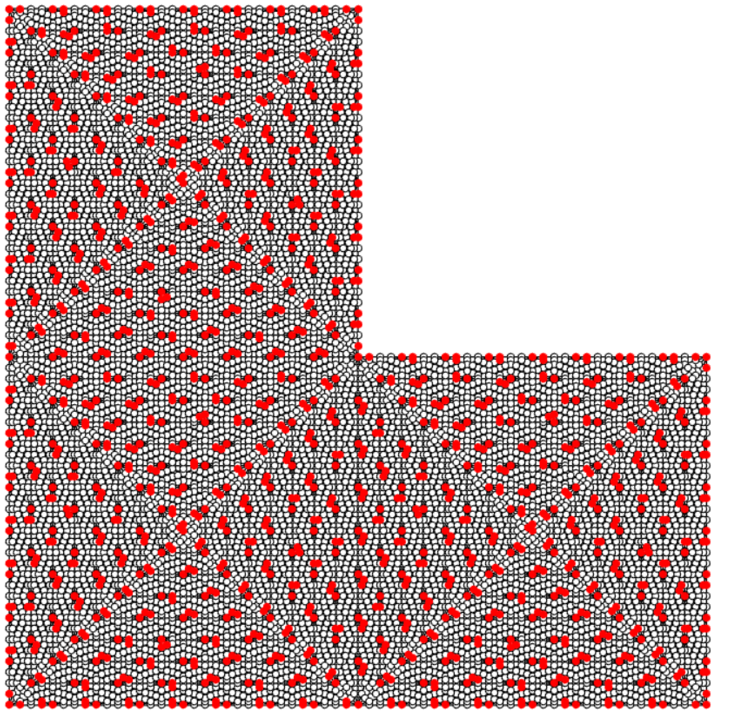}}\\
\subfigure[Results on $C^{0}$ Mesh]{\includegraphics[scale=0.35]{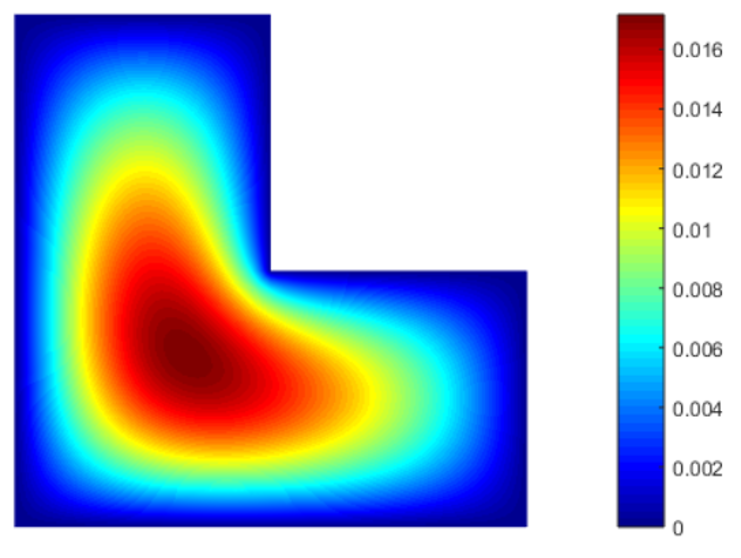}}\hspace{5mm}%
\subfigure[Results on $C^{1}$ Mesh]{\includegraphics[scale=0.35]{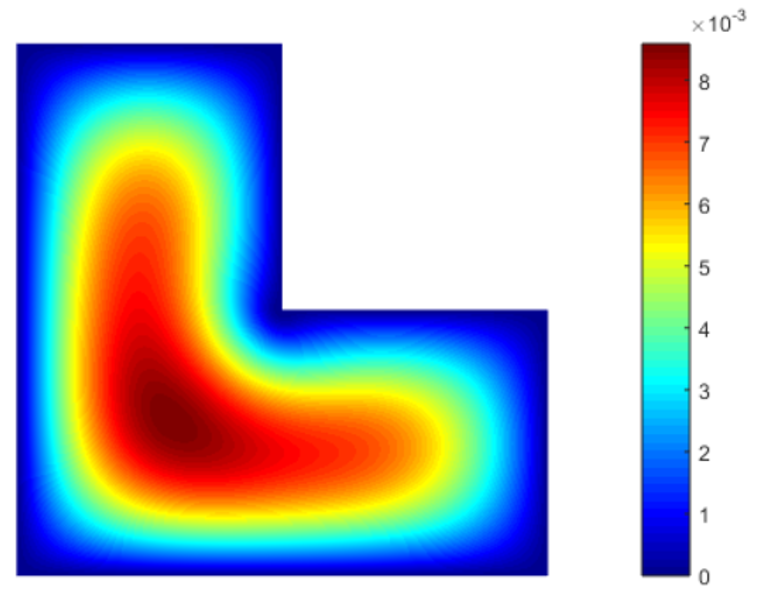}}
\caption{L-shaped non-convex domain, numerical results from Mixed FEM (left) and TIGA (right, created by PS) converged to different solutions. Results from Mixed FEM are 100\% larger than results from TIGA.}
\label{lshaped} 
\end{figure}

\begin{figure}
\centering
\subfigure[$C^{0}$ Mesh]{\includegraphics[scale=0.35]{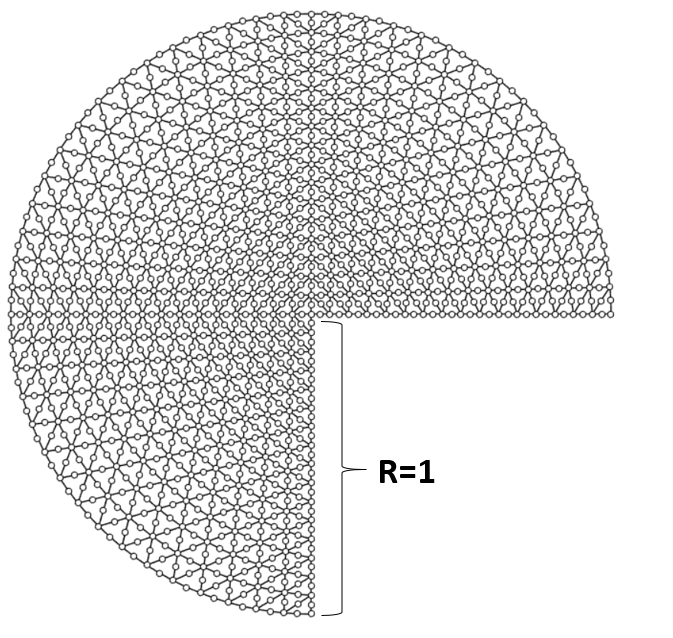}}\hspace{5mm}%
\subfigure[$C^{1}$ Mesh]{\includegraphics[scale=0.35]{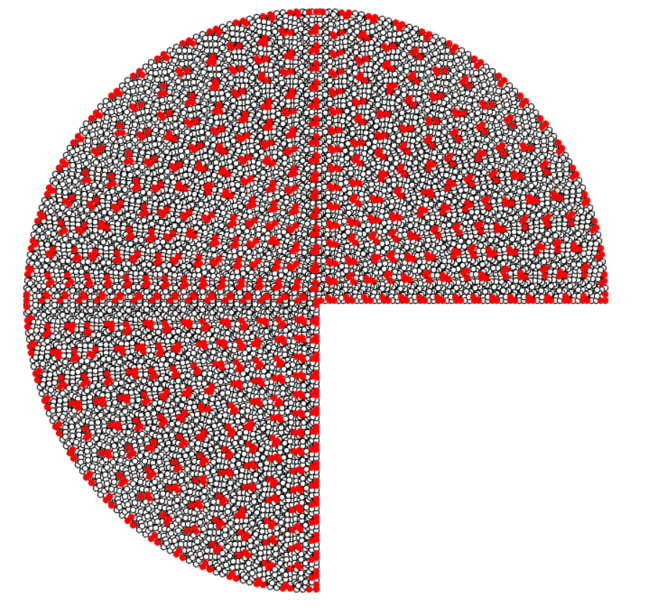}}\\
\subfigure[Results on $C^{0}$ Mesh]{\includegraphics[scale=0.22]{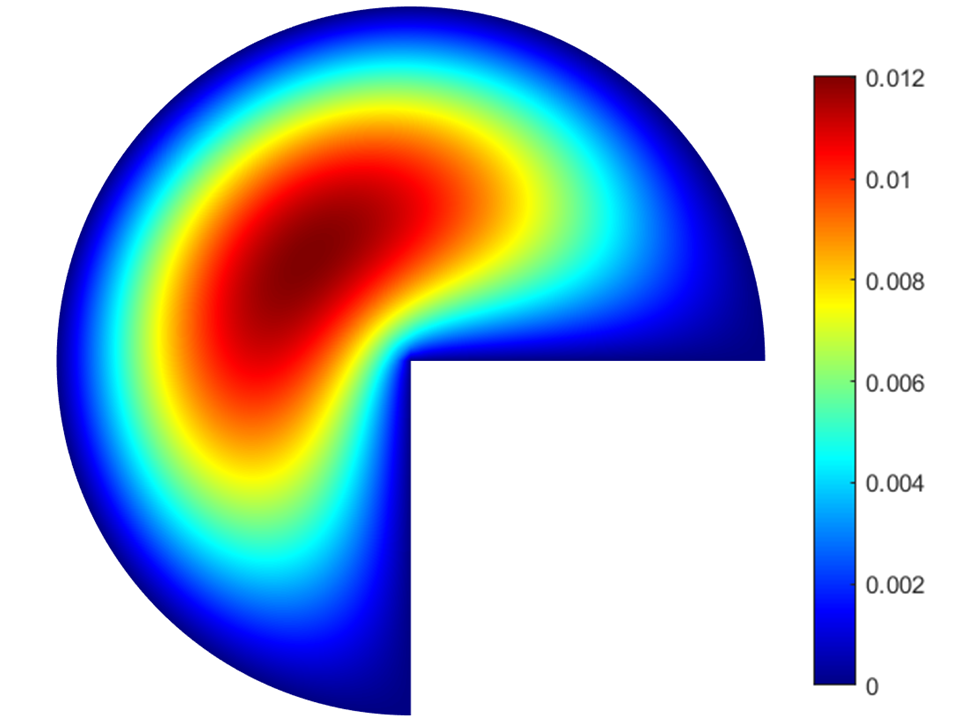}}\hspace{5mm}%
\subfigure[Results on $C^{1}$ Mesh]{\includegraphics[scale=0.218]{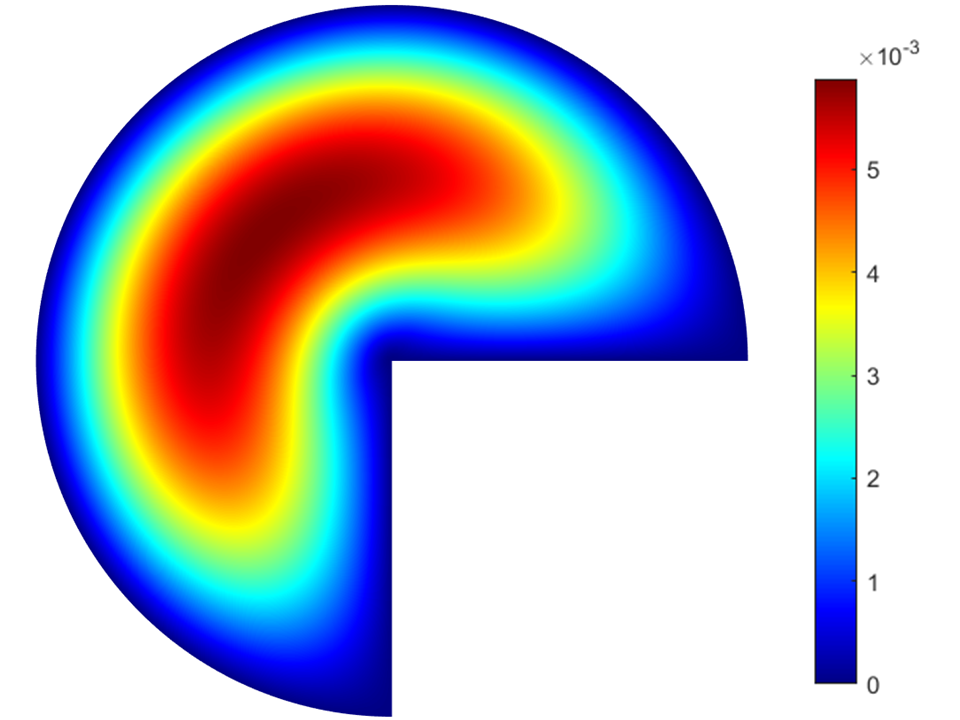}}
\caption{Circular L-shaped non-convex domain, numerical results from Mixed FEM (left) and TIGA (right, created by PS) converged to different solutions. Results from Mixed FEM are 100\% larger than results from TIGA.}
\label{clshaped} 
\end{figure}

\begin{figure}
\centering
\subfigure[$C^{0}$ Mesh]{\includegraphics[scale=0.35]{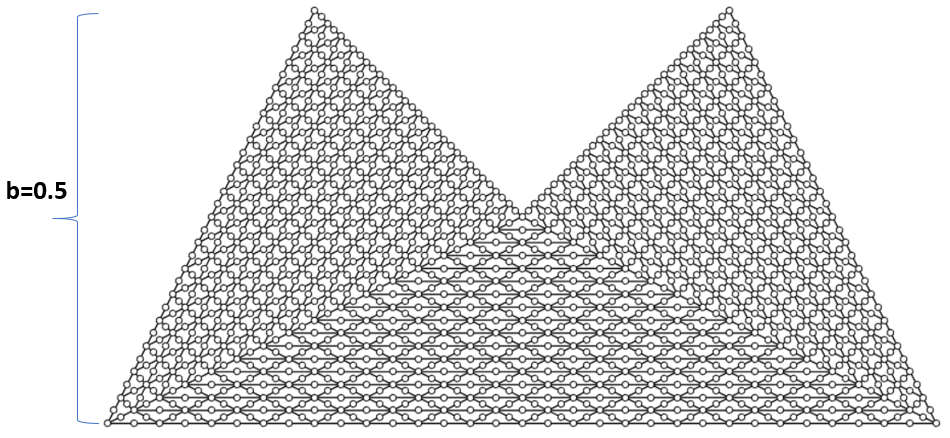}}\hspace{5mm}%
\subfigure[$C^{1}$ Mesh]{\includegraphics[scale=0.35]{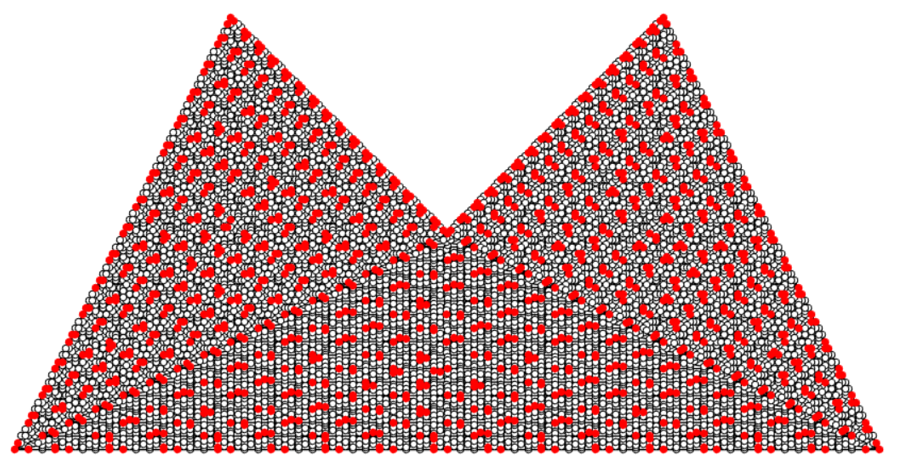}}\\
\subfigure[Results on $C^{0}$ Mesh]{\includegraphics[scale=0.4]{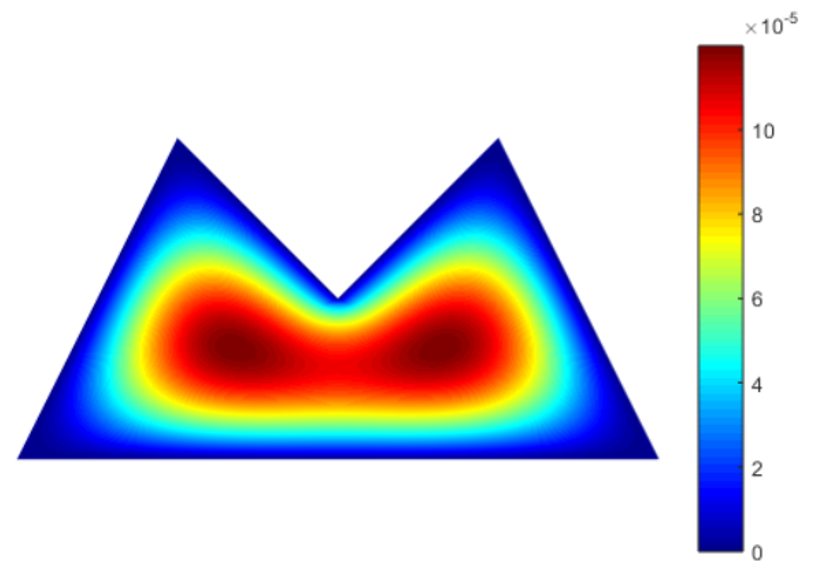}}\hspace{5mm}%
\subfigure[Results on $C^{1}$ Mesh]{\includegraphics[scale=0.4]{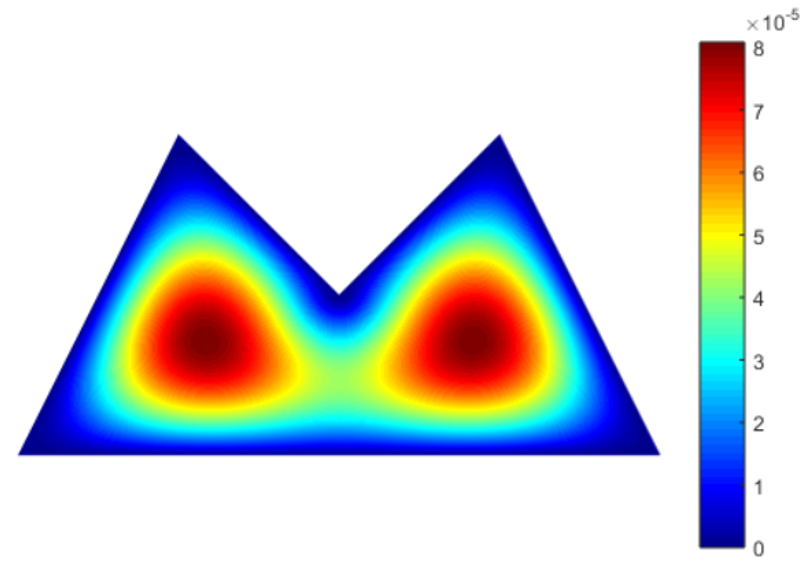}}
\caption{M-shaped non-convex domain,  numerical results from Mixed FEM (left) and TIGA (right, created by PS) converged to different solutions. Results from Mixed FEM are 50\% larger than results from TIGA.}
\label{Mshaped} 
\end{figure}

Figures \ref{lshaped}, \ref{clshaped} and \ref{Mshaped} compare the reliability of mentioned techniques. Each figure shows a geometry with re-entrant corner, $C^{0}$ mesh, $C^{1}$ mesh created by PS along with the results obtained from applying simply supported boundary condition and unity as the rhs function, $f$. In mixed FEM. all nodes are free; however, in $C^{1}$ TIGA, there are free nodes (red solid nodes) and dependent nodes (white nodes) in mesh.  

As it can be observed, non-convex domain with re-entrant corners, numerical results from mixed FEM(left) and TIGA (right) converged to different solutions; \hcb{ it is worth noting that we compare results after conducting mesh independence studies for each case.} In Figure \ref{lshaped}, maximum computed value of $u$ is $0.016$ and $0.008$ for mixed FEM (1601 nodes) and $C^{1}$ TIGA (1251 free nodes), respectively, i.e.  result from mixed FEM is 100\% larger than results from $C^{1}$ TIGA. Same relative difference is observed in Figure \ref{clshaped}(1275 free nodes in $C^{1}$ mesh and 1617 nodes in $C^{0}$ mesh). In Figure \ref{Mshaped}, maximum estimated $u$ is $12\times 10^{-5}$ and $8\times 10^{-5}$ for mixed FEM (1617 nodes) and $C^{1}$ TIGA (1275 free nodes), respectively. Moreover, in the same case, we can observe a slight difference in estimated distribution of $u$ over the domain. 

\section*{Conclusion}

In this study we implemented $C^{1}$ TIGA and $C^{0}$ TIGA (mixed FEM) for numerical solution of von Karman equations. Error analysis of the presented method demonstrate that we can obtain optimal convergence rate for $L^{2}$,$H^{1}$ and $H^{2}$ using $C^{1}$ TIGA, both quadratic and cubic elements.
	
Regarding efficiency in needed DOFs, TIGA is highly efficient. Mixed FEM is computationally expensive compared with $C^{1}$ TIGA for same accuracy, e.g. for cubic elements, mixed formulations could cost 490\% more DOFs for the same accuracy. This can also lead to convergence difficulties in nonlinear problems such as von Karman equations for geometrically nonlinear plate.

Furthermore, results of the current work demonstrate that when domain has re-entrant corner, mixed FEM can not converge to the solution obtained from TIGA ($C^{1}$ elements). The difference between two types of solutions is significant; same observation can be found in an article by \cite{Nazarov2007,Gerasimov2012a}. 

This study shows that TIGA is a promising tool in investigating more engineering problems for which consideration of high accuracy, efficiency and nonlinearity are of particular interest.

\section*{Acknowledgment}
The authors would like to acknowledge the financial support from ARO grant W911NF-17-1-0020.


\bibliographystyle{plain}
\bibliography{vonKarman_ref}

\end{document}